\documentclass[11pt]{amsart}

\usepackage{amssymb}
\usepackage{amsmath}

\usepackage{latexsym}
\usepackage{mathrsfs}          

\usepackage[all]{xy}           

\theoremstyle{plain} 
\newtheorem{theorem}{Theorem}[section]
\newtheorem{proposition}[theorem]{Proposition}

\newtheorem{corollary}[theorem]{Corollary}

\theoremstyle{definition} 
\newtheorem{definition}[theorem]{Definition}
\newtheorem{example}[theorem]{Example}

\theoremstyle{remark} 

\newcounter{ilist}
\newenvironment{ilist}{\begin{list}
                        {$(\mathrm{\roman{ilist}})$}
                        {\usecounter{ilist}\leftmargin0.9cm\labelwidth1cm}
                      }{\end{list}}

\newcounter{1list}
\newenvironment{1list}{\begin{list} 
                        {$( \mathrm{\arabic{1list}})$}
                        {\usecounter{1list} \itemsep0.8pt \topsep-1.7pt \leftmargin0.9cm
                          \labelwidth1cm}
                      }{\end{list}}

\newcommand{\bC}{\mathbb{C}}
\newcommand{\bF}{\mathbb{F}}

\newcommand{\bH}{\mathbb{H}}

\newcommand{\bP}{\mathbb{P}}
\newcommand{\bQ}{\mathbb{Q}}

\newcommand{\bV}{\mathbb{V}}
\newcommand{\bZ}{\mathbb{Z}}

\newcommand{\cA}{\mathcal{A}}

\newcommand{\cM}{\mathcal{M}}

\newcommand{\cO}{\mathcal{O}}

\newcommand{\C}{\bC}  
\newcommand{\F}{\bF}  
\newcommand{\Q}{\bQ}  
\newcommand{\Z}{\bZ}  

\newcommand{\Wlog}{Without loss of generality}

\newcommand{\pp}{principally polarized}
\newcommand{\ppav}{principally polarized Abelian variety}
\newcommand{\ppavs}{principally polarized Abelian varieties}
\newcommand{\ppass}{principally polarized Abelian schemes}

\newcommand{\Fq}{{\F_q}}

\newcommand{\Fqn}{\F_{q^n}}

\newcommand{\omegacy}{\omega_{C/Y}}

\newcommand{\rar}[1]{\stackrel{#1}{\rightarrow}}
\newcommand{\grar}[1]{\stackrel{#1}{\longrightarrow}}
\newcommand{\raro}{{\rightarrow}}

\newcommand{\Kbar}{\overline{K}}
\newcommand{\wh}[1]{\widehat{#1}}

\newcommand{\GK}{{G_{\bar{K}/K}}}
\newcommand{\Pone}{{\bP^1_\C}}
\newcommand{\Ag}{\cA_{g}}
\newcommand{\Mg}{\cM_g}

\newcommand{\isom}{\simeq}

\newcommand{\Fr}[1]{F_#1}
\newcommand{\Hl}[1]{H^1(\bar{#1}, \Q_\ell)}
\newcommand{\TrF}[2]{\mathrm{Tr}\big(\Fr{#1}^{#2}\, \big|\, H^1(\bar{#1}, \Q_\ell)\big)}
\newcommand{\Gal}[1]{Gal\big(\overline{#1}/#1\big)}

\newcommand{\Spec}{\mathrm{Spec}\,}
\newcommand{\End}{\mathrm{End}}

\newcommand{\Image}{\mathrm{Im}}
\newcommand{\Ker}{\mathrm{Ker}}

\newcommand{\GL}{\mathrm{GL}}
\newcommand{\SL}{\mathrm{SL}}

\begin{document}


\title[On Shimura curves in the Schottky locus]{On Shimura curves in the Schottky locus}
\author[S. Kukulies]{Stefan Kukulies}
\address{Universit\"at Duisburg-Essen, Mathematik, 45117 Essen, Germany}
\email{Stefan.Kukulies@uni-due.de}

\maketitle


\begin{abstract}
  We show that a given rational Shimura curve $Y$ with strictly maximal Higgs field in the moduli
  space of $g$-dimensional Abelian varieties does not generically intersect the Schottky locus for
  large $g$.

  We achieve this by using a result of Viehweg and Zuo which says that if $Y$ parameterizes a family
  of curves of genus $g$, then the corresponding family of Jacobians is $Y$-isogenous to the
  $g$-fold product of a modular family of elliptic curves. After reducing the situation from the
  field of complex numbers to a finite field, we will see, combining the Weil and Sato-Tate
  conjectures, that this is impossible for large $g$.  
\end{abstract}


\tableofcontents


\section{Introduction}

Let $U$ be a rational Shimura curve so that $U$ is an open subset in $\Pone$. Let $S=\Pone - U$ be
the ``bad locus''. We prove the following theorem.\pagebreak[3]

\begin{theorem}[Shimura curves in the Schottky locus]\label{theorem1}
  Given an integer $s\geq 0$, there is a natural number $B=B(s)$, depending only on $s$, such that a
  rational Shimura curve, whose Higgs field is strictly maximal and whose bad locus $S$ contains at
  most $s$ points, cannot lie in the closure of the Schottky locus $\Mg$ for $g>B$.
\end{theorem}\pagebreak[3]

The Schottky locus is the image of the moduli space $\Mg$ of curves of genus $g$ in $\Ag$ - the
moduli space of principally polarized $g$-dimensional Abelian varieties (with a suitable level
structure). We say that a Shimura curve in $\Ag$ lies in the closure of the Schottky locus if it
generically intersects the image of $\Mg$ in $\Ag$. So, in particular, it should not
lie entirely in the boundary of the closure of $\Mg$ in $\Ag$.\pagebreak[3]

By Shimura variety we mean a Shimura variety of Hodge type which is an \'etale covering of a certain
moduli space of Abelian varieties with prescribed Mumford-Tate group and a suitable level structure
as defined in \cite{MumfordShimura}. A Shimura curve is a one-dimensional Shimura
variety.\pagebreak[3]

Let $f:A\raro Y$ be a semistable family of Abelian varieties over a complex projective curve $Y$, $U=Y-S$ the
smooth locus and $V=f^{-1}(U)$ so that $f:V\raro U$ is an Abelian scheme. Consider the Higgs bundle
$(E,\theta)$ given by taking the graded sheaf of the Deligne extension of $R^1f_*\C_V\otimes\cO_U$
where $R^1f_*\C_V$ is the weight $1$ variation of Hodge structures. We have a decomposition
$E=F\oplus N$ into an ample part $F$ and a flat part $N$. Following \cite{VZ-strictMaxHiggs} we say
that the Higgs field is {\it maximal} if
$$ \theta^{1,0}: F^{1,0}\longrightarrow F^{0,1}\otimes \Omega^1_Y(\log S) $$
is an isomorphism, and that the Higgs field is {\it strictly maximal} if additionally
$N=0$.\pagebreak[3]

Viehweg and Zuo showed in \cite{VZ-CharShCurves} that if each irreducible and non-unitary
sub-variation $\bV$ of Hodge structures in $R^1f_*\C_V$ has a strictly maximal Higgs field, then
there is an \'etale covering $U'\raro U$ such that $U'$ is a Shimura curve and $f':V'\raro U'$ is
the corresponding universal family. Moreover, M\"oller showed in \cite{Moeller2} that the converse
also holds. Hence we have a characterization of Shimura curves by the maximality of the Higgs field
of the corresponding universal family.\pagebreak[3]

Combining the results of Viehweg and Zuo \cite{VZ-NumericalBounds} which say that a Shimura curve
$U$ in $\Mg$ has to be non-compact with the techniques of M\"oller shows that $U$ has also to be a
Teichm\"uller curve. Then from \cite{Moeller2} it follows that there are no such curves in $\Mg$
unless $g=3$. See also the discussion in \cite{MVZ}.\pagebreak[3]

Observe that this result deals with the occurrence of Shimura curves in $\Mg$ rather than its closure
in $\Ag$. So it does not answer the question if there are Shimura curves in the closure of the
Schottky locus.\pagebreak[3]

We remark that the conjecture of Andr\'e-Oort, saying that a Shimura variety is characterized by
having a dense set of CM-points, and the conjecture of Coleman, saying that there are only finitely
many CM-points in $\Mg$ for high genus $g$, suggest that there are no Shimura varieties in the
closure of the Schottky locus for $g$ sufficiently large.\pagebreak[3]

In \cite{Hain}, Hain studied families of compact Jacobians over locally symmetric domains $U$
satisfying an additional technical condition. Based on his methods, de Jong and Zhang
\cite{deJongZhang} were able to exclude certain types of higher-dimensional Shimura
varieties.\pagebreak[3]

Returning to Shimura curves, we will prove Theorem (\ref{theorem1}) as follows. Let $C\raro Y$ be a
family of complex curves whose family of Jacobians $J\raro Y$ has a strictly maximal Higgs field. If
$Y=\Pone$, then a further result of Viehweg and Zuo from \cite{VZ-CharShCurves} says that $J\raro Y$
is $Y$-isogenous to the $g$-fold product $E\times_Y\cdots\times_Y E$ of a modular family of elliptic
curves $E\raro Y$. So Theorem (\ref{theorem1}) will follow from the following theorem, which holds
for an arbitrary base curve $Y$.\pagebreak[3]

\begin{theorem}[Bound for the genus]\label{theorem2}
  Let $C\raro Y$ be a family of curves of genus $g$ whose Jacobian $J\raro Y$ is $Y$-isogenous to
  the $g$-fold product of a non-isotrivial family of elliptic curves $E\raro Y$ which can be defined
  over a number field. Then the genus $g$ is bounded, i.\,e. there is a number $B=B(E/Y)$,
  depending only on $E\raro Y$, such that $g$ is smaller than $d$.
\end{theorem}\pagebreak[3]

Mind that modular families of elliptic curves can be defined over number fields. We will prove
Theorem (\ref{theorem2}) by reducing the situation from $\C$ to a number field $F$. Then we reduce to
a finite field by selecting a suitable finite  prime of $F$. Finally, we prove that the genus $g$ of
the fibers of a family of curves $C\raro Y$, defined over a finite field and whose family of
Jacobians is isogenous to the $g$-fold product of a family of elliptic curves, is
bounded.\pagebreak[3]

We achieve this by counting the number of singularities $\delta$ in the fibers of $C\raro
Y$. Combining the Weil conjectures for the fibers with the Sato-Tate conjecture about the
distribution of Frobenius traces in a family of elliptic curves, we will get a lower bound for
$\delta$. On the other hand, the geometry of the total space $C$ of $C\raro Y$ will give an upper
bound. For large $g$, the lower bound will exceed the upper bound. Thus, the genus has to be
bounded.\pagebreak[3]\newline

I would like to thank my adviser Eckart Viehweg for his continuous support and encouragement while
writing this paper. I would also like to thank Gebhard B\"ockle, Martin M\"oller and Pooji Singla
for pointing out several mistakes. This work was financially supported by the DFG.\pagebreak[3]


\section{Shimura curves and maximal Higgs fields}

Let $Y$ be an irreducible smooth projective curve over the field of complex numbers $\bC$ and
let $f:A\raro Y$ be a semistable family of $g$-dimensional Abelian varieties, i.\,e. $f:A\raro Y$ is
a flat, projective $\bC$-morphism whose generic fiber is an Abelian variety of dimension $g$. Let
$U\subset Y$ be the smooth locus of $A\raro Y$, i.\,e. the restriction of $A\raro Y$ to $U$ is an
Abelian scheme $A_0\raro U$ while the fibers over the set $S=Y - U$ are all singular. Consider the
weight $1$ variation of Hodge structures $R^1f_*\Z_{A_0}$ and let $F$ be the non-flat part of the
Higgs bundle $(E,\theta)$ given by taking the graded sheaf of the Deligne extension of
$R^1f_*\Z_{A_0}\otimes\cO_U$ to $Y$ which carries a Hodge filtration. Then the Arakelov inequality
for families of Abelian varieties \cite{JostZuo} says that 
$$ 0< 2\cdot \deg(F^{1,0}) \leq g_0\cdot \big(2q-2+\#S\big) $$
where $q$ denotes the genus of the base curve $Y$ and $g_0$ is the rank of $F^{1,0}$. We say that
the family of Abelian varieties $A\raro Y$ {\it reaches the Arakelov bound } if the above inequality
becomes an equality. Viehweg and Zuo showed in \cite{VZ-CharShCurves} that this property is
equivalent to the {\it maximality} of the Higgs field for $F$, i.\,e. the map $\theta|_{F^{1,0}}:
F^{1,0} \raro F^{0,1}\otimes \Omega^1_Y(\log S)$ is an isomorphism. Moreover, the Higgs field is
called {\it strictly maximal} if in addition the Higgs bundle has no flat part.\pagebreak[3]

Assume that $V\raro U$ is a Shimura curve, i.\,e. $U$ is an \'etale covering of a certain moduli
space of Abelian varieties with prescribed Mumford-Tate group and a suitable level structure, and
$V\raro U$ is the corresponding universal family, see \cite{Moeller2}. Then M\"oller showed that
$U\raro V$ has a maximal Higgs field.\pagebreak[3]

\begin{theorem}[Shimura implies maximal Higgs]
  If $V\raro U$ is the universal family over a Shimura curve, then its Higgs field is maximal.
\end{theorem}\pagebreak[3]

\begin{proof}
  See \cite[Thm.0.9]{MVZ} or \cite[Thm.1.2]{Moeller2}.
\end{proof}\pagebreak[3]

The converse was shown by Viehweg and Zuo in \cite{VZ-CharShCurves}, see also \cite{MVZ}. So we have
a Characterization of Shimura curves by its corresponding Higgs field. Moreover, Viehweg and Zuo
showed that after an \'etale extension the family $A\raro Y$ decomposes in the following
way.\pagebreak[3]

\begin{theorem}[Decomposition Theorem]\label{thm2.2}
  If $V\raro U$ has a maximal Higgs field and $S\neq\emptyset$, the there is an \'etale covering
  $Y'\raro Y$ such that the pull-back family $A'\raro Y'$ is $Y'$-isogenous to a product
  $$ E\times_{Y'}\cdots\times_{Y'} E\times_\bC B $$
  where $B/\bC$ is an Abelian variety of dimension $g-g_0$ and $E\raro Y'$ is a modular family of
  elliptic curves.
\end{theorem}\pagebreak[3]

\begin{proof}
  See \cite[Thm.0.2]{VZ-CharShCurves} or \cite[Cor.0.10]{MVZ}.
\end{proof}\pagebreak[3]

Modular means that the smooth locus $U'$ of $E\raro Y'$ is the
quotient $\Gamma\backslash\bH$ of the upper half-plane $\bH$ by a subgroup $\Gamma\subset\SL_2(\Z)$
of finite index and $E\raro Y'$ is over $U'$ the quotient of $\bH\times\C$ by the semi-direct
product of $\Gamma$ and $\Z^2$.\pagebreak[3]

If the Higgs field is strictly maximal so that $g_0=g$, then there is no constant part. Hence,
$A'\raro Y'$ is $Y'$-isogenous to a product $E\times_{Y'}\ldots\times_{Y'}E$ with $E\raro Y'$
modular.\pagebreak[3]

Now let $C\raro Y$ be a semistable family of curves of genus $g$ with $J\raro Y$ its corresponding
family of Jacobians. Let $V\raro U$ be the the smooth part of $J\raro Y$ and $S=Y - U$ the
bad locus - as before with $A\raro Y$ instead of $J\raro Y$. Assume that $V\raro U$ is a Shimura
curve so that its Higgs field is maximal. If $Y=\Pone$, then the Arakelov inequality above tells us
that $S\neq\emptyset$. If we further assume that the Higgs field of $V\raro U$ is not only maximal
but strictly maximal, then the Jacobian has the following decomposition.\pagebreak[3]

\begin{corollary}[Structure of the Jacobian]\label{corollarySplit}
  Let $C\raro Y$ be a semistable family of curves of genus $g$ whose Jacobian $J\raro Y$ has a
  strictly maximal Higgs field. If $Y=\Pone$, then $J\raro Y$ is $Y$-isogenous to the $g$-fold
  product $E\times_Y\cdots\times_Y E$ of a modular family of elliptic curves $E\raro Y$.
\end{corollary}\pagebreak[3]

\begin{proof}
  Because of the Arakelov inequality we have $S\neq\emptyset$. So we may apply Theorem
  (\ref{thm2.2}). Since $Y=\Pone$ has no other \'etale coverings than automorphisms, the splitting
  takes place over $Y$.
\end{proof}\pagebreak[3]

We will show that for a fixed modular family $E\raro Y$ there cannot exist families of curves
$C\raro Y$ of arbitrary large genus whose Jacobian is a $g$-fold product of $E\raro Y$.\pagebreak[3]


\section{Reduction to number fields}

We want to show that a family of curves $C\raro Y$ whose Jacobian $J\raro Y$ is $Y$-isogenous to the
$g$-fold product of a modular family of elliptic curves $E\raro Y$ is (after a base change) defined
over a number field $F$ which depends only on $E\raro Y$ and not on $g$.\pagebreak[3]

We say that a family of curves or group schemes $X\raro Y$ is {\it defined over a number field $F$}
if there is a curve $Y_0/F$ and a family of curves or group schemes $X_0\raro Y_0$ over $F$ such that
$Y_0\times_F\bC \isom Y$ over $\bC$ and the pull-back family $X_0\times_F\bC\raro
Y_0\times_F\bC\isom Y$ is $Y$-birational to $X\raro Y$, i.\,e. the families have isomorphic generic
fibers.\pagebreak[3]

We say that a $Y$-morphism $f$ between two families $X\raro Y$ and $Z\raro Y$ is {\it defined over a
  number field $F$} if $X\raro Y$ and $Z\raro Y$ are defined over $F$ and there is a $Y_0$-morphism
$f_0:X_0\raro Y_0$ which coincides generically with the $Y$-morphism $f:X\raro Y$ after the base
change $\Spec\bC \raro\Spec F$.\pagebreak[3]

As we will see, the reason why $C\raro Y$ descends to a number field is that modular families of
elliptic curves $E\raro Y$ are defined over number fields. We start by describing the torsion
structure of families of elliptic curves $E\raro Y$ via Galois representations.\pagebreak[3]

Let $K:=\bC(Y)$ be the function field of $Y$ and $K_v$ be the $v$-adic completion of $K$ where $v$
denotes a normalized discrete valuation of $K$ induced by some point $y\in Y(\bC)$. Let $G_v =
\Gal{K_v}$ be the absolute Galois group of $K_v$. By $j_E$ we denote the $j$-invariant of an
elliptic curve.\pagebreak[3]

\begin{proposition}[Galois action on torsion of Tate curves]\label{torsionOnTate}
  Let $K_v$ be a $v$-adic complete field with residue field $\C$  and absolute Galois group $G_v$,
  and let $E/K_v$ be a Tate curve. Then for any prime power $\ell^n$, we can find a basis
  $(P_1,P_2)$ of $E[\ell^n]\big(\overline{K_v}\big)$ such that for any integer $n'$ with $\ell^{n'+1}\nmid
  v(j_E)$, there is an element $\sigma\in G_v$ which acts on $E[\ell^n]\big(\overline{K_v}\big)$ with respect
  to the basis $(P_1,P_2)$ like 
  $$\left(\begin{array}{ccc} 
      1 & \ell^{n'}\\
      0 & 1
    \end{array}\right)\in\GL_2(\Z/\ell^n\Z).
  $$
  In particular, for almost all prime powers $\ell^n$ there is a transvection, i.\,e. $n'=0$.

  Furthermore, the basis $(P_1,P_2)$ can be chosen as follows: for $P_1$ we may take any
  $\ell^n$-torsion point which specializes into the connected component of one, while for $P_2$ we
  may take any other point such that $(P_1,P_2)$ forms a basis.
\end{proposition}\pagebreak[3]

\begin{proof}
  Mimic the proof of \cite[V.6.1]{Silverman2} using Tate's $v$-adic Uniformization Theorem.
\end{proof}\pagebreak[3]

From this we can conclude that the image of the action of the absolute Galois group $G=\Gal{K}$
acting on $N$-torsion points $E[N]\big(\Kbar\big)$ of $E\raro Y$ is huge.\pagebreak[3]

\begin{corollary}[Galois action on torsion points of $E$]\label{GaloisOnTorsion}
  Let $E\raro Y$ be a non-isotrivial family of elliptic curves and $K=\bC(Y)$ the function field of
  $Y$. Then for any prime number $\ell$ there is a non-negative integer $n(\ell)$ such that for
  all prime powers $\ell^n$, there are elements $\sigma$ and $\sigma'$ of $G$ which act like
  $$\left(\begin{array}{ccc} 
      1 & \ell^{n(\ell)}\\
      0 & 1
    \end{array}\right)\quad \textrm{and}\quad
    \left(\begin{array}{ccc} 
      1 & 0\\
      \ell^{n(\ell)} & 1
    \end{array}\right)
  $$
  on $E[\ell^n]\big(\overline{K}\big)$ with respect to a suitable basis. Moreover, for almost all
  $\ell$, we may choose $n(\ell)=0$.
\end{corollary}\pagebreak[3]

\begin{proof}
  We argue as in \cite{Igusa3}.
  We may assume that $E\raro Y$ has everywhere semistable reduction and a full
  level-$\ell^n$-structure, i.\,e. there is an isomorphism of $Y$-group schemes
  $(\Z/\ell^n\Z)^2_Y\raro E[\ell^n]$. This can always be achieved after a finite base
  change.\pagebreak[3]

  Choosing a point $y\in Y(\C)$ such that $E\raro Y$ has bad reduction in $y$, we find by
  Proposition (\ref{torsionOnTate}) a basis $(P_1, P_2)$ of $E[\ell^n](\Kbar)$ such that there is an
  element  $\sigma \in G$ which acts with respect to $(P_1, P_2)$ like 
  $$\left(\begin{array}{ccc} 
      1 & \ell^{n'}\\
      0 & 1
    \end{array}\right)\in\GL_2(\Z/\ell^n\Z)
  $$
  for $n'$ with $\ell^{n'+1}\nmid v(j_E)$ where $v$ is the valuation at $y$. Furthermore, $P_1$
  specializes into the connected component of one while $P_2$ does not. So we may find another point
  $y' \in Y(\C)$ such that $P_2$ will specialize into the connected component of one, since having a
  full level-$\ell^n$-structure $E\raro Y$ is the pull-back of the universal elliptic curve
  $E(\ell^n)\raro X(\ell^n)$ parameterizing full level-$\ell^n$-structures (we may assume that
  $\ell^n>2$ because if the statement is true for $\ell^n$, it is also true for $\ell^{n-1}$). So
  using again (\ref{torsionOnTate}) we will find an element $\sigma'\in G$ which acts with respect to
  the basis $(P_2,P_1)$ like 
  $$\left(\begin{array}{ccc} 
      1        & \ell^{n''}\\
      0        & 1
    \end{array}\right)\in\GL_2(\Z/\ell^n\Z)
  $$
  for $n''$ with $\ell^{n''+1}\nmid v'(j_E)$ where $v'$ is the valuation in $y'$. Of course upper
  triangle matrices with respect to $(P_2,P_1)$ will be lower triangle matrices with respect to
  $(P_1,P_2)$. So choosing $n(\ell)$ such that $\ell^{n(\ell)+1} \nmid v(j_E)$ and
  $\ell^{n(\ell)+1} \nmid v'(j_E)$, we find two elements $\sigma$ and $\sigma'$ of $\GK$ which act
  with respect to $(P_1,P_2)$ like
  $$\left(\begin{array}{ccc} 
      1 & \ell^{n(\ell)}\\
      0 & 1
    \end{array}\right)\quad \textrm{and}\quad
    \left(\begin{array}{ccc} 
      1 & 0\\
      \ell^{n(\ell)} & 1
    \end{array}\right).
  $$
  In particular, $n(\ell)$ does only depend on $\ell$ or $\ell^n$ and for almost all $\ell$,
  we may choose $n(\ell)=0$ because $\ell\nmid v(j_E)$ and $\ell\nmid v'(j_E)$.
\end{proof}\pagebreak[3]

We can conclude the shape of endomorphisms on torsion groups of $E\raro Y$.\pagebreak[3]

\begin{corollary}[Endomorphisms of torsion groups $E\lbrack N\rbrack$]\label{cor12}
  Let $E\raro Y$ be a non-isotrivial family of elliptic curves. Then for any prime number
  $\ell$, there is a non-negative integer $n(\ell)$ such that for all prime powers $\ell^n$,
  every   $Y$-endomorphisms  of $E[\ell^{n}]$ is the sum of a multiplication-by-$m$ map and a
  composition of the multiplication-by-$\ell^{n-n(\ell)}$ with an endomorphism of
  $E[\ell^{n(\ell)}]$.

  Moreover, for almost all $\ell$ we may choose $n(\ell)=0$ so that the endomorphisms
  $\End_Y\big(E[\ell^{n}]\big)$ consist only of multiplication-by-$m$ maps.
\end{corollary}\pagebreak[3]

\begin{proof}
  Any endomorphism of $E[\ell^n]$ is invariant under the action of Galois so that it has to lie
  in the center of the action of Galois on $E[\ell^n]$. In particular, an endomorphism must commute
  with the two matrices from Corollary (\ref{GaloisOnTorsion}). An elementary matrix calculation
  shows that any such endomorphism has to be represented by a matrix of the form
  $$ A = m\cdot I + \ell^{n-n(\ell)}\cdot M \in \mathrm{M}_\ell(\Z/\ell^n\Z) $$
  where $m$ is an integer, $I$ is the identity matrix and $M$ is some other
  matrix. The interpretation is the following. $m\cdot I$ corresponds to the multiplica\-tion-by-$m$
  map, while $\ell^{n-n(\ell)}\cdot M$ is the composition of the
  multiplication-by-$\ell^{n-n(\ell)}$ map $E[\ell^n]\raro E[\ell^{n(\ell)}]$ with an endomorphism
  of $E[\ell^{n(\ell)}]$.\pagebreak[3]

  Moreover, Corollary (\ref{GaloisOnTorsion}) says that for almost all $\ell$ we have $n(\ell)=0$ so
  that $\End_Y\big(E[\ell^{n(\ell)}]\big) = \{0\}$.
\end{proof}\pagebreak[3]

This tells us that if the family $E\raro Y$ is defined over some number field, then the same is true
for the endomorphisms of $E[N]$.\pagebreak[3]

\begin{corollary}[Endomorphisms of torsion groups descend]\label{cor13}
  Let $E\raro Y$ be a non-isotrivial family of elliptic curves defined over some number field
  $F$. Then there is a finite extension $F'$ of $F$ such that for all natural 
  numbers $N$ the $Y$-endomorphisms of $E[N]$ are defined over $F'$. In particular, the field $F'$
  depends only on $E\raro Y$.
\end{corollary}\pagebreak[3]

\begin{proof}
  It is enough to consider prime powers $\ell^n$. By Corollary (\ref{cor12}) for almost all $\ell$
  the endomorphisms $\End_Y\big(E[\ell^{n(\ell)}]\big)$ consist only of multiplication-by-$m$ maps
  which are clearly defined over $F$.\pagebreak[3]

  For the finitely many remaining $\ell$, we also have to consider endomorphisms of
  $E[\ell^{n(\ell)}]$. Since these are finite in number, they will be defined over some finite
  extension $F'$ of $F$ depending only on $E\raro Y$.
\end{proof}\pagebreak[3]

In particular, any family of Abelian varieties $A\raro Y$ which is isogenous to a $g$-fold product of
$E\raro Y$ can be defined over the number field $F'$.\pagebreak[3]

\begin{proposition}[Isogenies and Abelian varieties descend]\label{pro14}
  Let $E\raro Y$ a non-isotrivial family of elliptic curves defined over some number field $F$. Then
  there is a finite covering $Y'\raro Y$ and a finite extension $F'$ of $F$ such that for every
  $Y$-isogeny $h$ from any $g$-fold product of $E\raro Y$ to any family of Abelian varieties $A\raro
  Y$, the $Y'$-isogeny $h'=h\times_Y id_{Y'}$ is defined over $F'$. In particular, $Y'$ and $F'$
  depend only on $E\raro Y$ and not on $g$.
\end{proposition}\pagebreak[3]

\begin{proof}
  Let $H$ be the kernel of $h$. For a suitable number $N$, the group scheme $H$ is
  contained in $E[N]\times_Y\cdots\times_Y E[N]$. After a suitable base change $Y'\raro Y$, which
  may be chosen independently from $H$ because of the finiteness result Corollary (\ref{cor12}), we
  can describe $H'$ as the kernel of an endo\-morphism of $E'[N]\times_{Y'}\cdots\times_{Y'}
  E'[N]$. Corollary (\ref{cor13}) shows that $H'$ will be defined over $F'$, i.\,e. $H'$ is the
  pull-back of a subgroup scheme $H'_0$ of $E'_0[N]\times_{Y'_0}\cdots\times_{Y'_0} E'_0[N]$ with
  respect to the base change $\Spec\bC\raro\Spec F'$.\pagebreak[3]

  Now let $A'_0\raro Y'_0$ be the quotient of the $g$-fold product of $E'_0\raro Y'_0$ by $H'_0$ and
  let $h'_0$ be the quotient map. Clearly, $h'$ and $A'\raro Y'$ coincide generically with the
  pull-backs of $h'_0$ and $A'_0$ under the base change $\Spec\bC\raro\Spec F'$. Also Corollary
  (\ref{cor13}) says that $F'$ depends only on the family $E'\raro Y'$ and, therefore, on $E\raro
  Y$.
\end{proof}\pagebreak[3]

Moreover, we can say something about the isogenies which will occur.\pagebreak[3]

\begin{proposition}[Divisibility of the degree of isogenies]\label{pro15}
  Let $E\raro Y$ be a non-isotrivial family of elliptic curves. Then there is a finite set of
  primes $S=S(E/Y)$, depending only on $E\raro Y$, such that for every family of Abelian varieties
  $A\raro Y$, which is $Y$-isogenous to a $g$-fold product of $E\raro Y$, there is a $Y$-isogeny
  between the $g$-fold product of $E\raro Y$ and $A\raro Y$ whose degree has only prime divisors
  contained in $S$. In particular, the set $S$ does not depend on $g$.
\end{proposition}\pagebreak[3]

\begin{proof}
  Let $S=S(E/Y)$ be the set of primes $\ell$ such that the integers $n(\ell)$ from Corollary
  (\ref{GaloisOnTorsion}) are non-zero. Let $h: E\times_Y\cdots\times_Y E\raro A$ be an isogeny and
  $H\subset \Ker(h)$ a non-trivial simple subgroup scheme so that the order of $G$ is some prime
  power $\ell^n$.\pagebreak[3]

  Let $\ell\notin S$. If $g=1$, then $H\isom E[\ell]$ since $E[\ell]$ is irreducible by Corollary
  (\ref{cor12}). Thus $h:E\raro A$ factorizes through $E\isom E/H\raro A$. Proceeding like this, we
  can find an isogeny prime to all $\ell\notin S$.\pagebreak[3]

  If $g>1$, then using Corollary (\ref{cor12}) we see that $H$ maps into a factor of
  $E\times_Y\cdots\times_Y E$ by a multiplication-by-$m$ map. So after applying a suitable
  automorphism of $E\times_Y\cdots\times_Y E$, the subgroup $H$ lies in a $g-1$-dimensional factor
  of the product. By induction we see that $(E\times_Y\cdots\times_Y E)/H$ is isomorphic to
  $E\times_Y\cdots\times_Y E$.
\end{proof}\pagebreak[3]

So far, we have seen that a family of Abelian varieties which is isogenous to a product of a modular
family of elliptic curves can be defined over a number field. Now we want to have this result for a
family of curves whose Jacobian is isogenous to a product of a modular family of elliptic
curves.\pagebreak[3]

\begin{proposition}[Curves descend]\label{uniformDescend}
  Let $E\raro Y$ be a non-isotrivial family of elliptic curves defined over
  a number field $F$. Then there is a finite field extension $F'$ of $F$ and a curve $Y'$
  covering $Y$, both depending only on $E\raro Y$, such that for any family of curves
  $C\raro Y$, whose Jacobian $J\raro Y$ is $Y$-isogenous to a $g$-fold product of
  $E\raro Y$, there is finite covering $Y''\raro Y'$ of degree at most $2$, such that
  $C\times_Y Y''\raro Y''$ is defined over $F'$. In particular, $Y'$ and the degree of
  $Y''\raro Y'$ do not depend on $g$.
\end{proposition}\pagebreak[3]

\begin{proof}
  Let $S=S(E/Y)$ be the set of primes from Proposition (\ref{pro15}) about the divisibility of the
  degree of isogenies. Let $N\geq 3$ be an integer which is not divisible by any prime in $S$. (The
  choice of $N$ will later ensure that the Jacobian is equipped with a level-$N$-structure, see
  below.)\pagebreak[3]

  After a finite extension $Y'\raro Y$, we may assume that $E'\raro Y'$ is equipped with a
  level-$N$-structure $\alpha'$. By further extending $Y'\raro Y$ we may also assume that the
  conclusion of Proposition (\ref{pro14}) about descending isogenies and Abelian varieties
  holds. The extension $Y'\raro Y$ depends only on $E\raro Y$.\pagebreak[3]

  Now let $C\raro Y$ be a curve whose Jacobian $J\raro Y$ is $Y$-isogenous to a $g$-fold
  product of $E\raro Y$. We may assume by Proposition (\ref{pro15}) that there is a $Y$-isogeny $h$
  from $E\times_Y\cdots\times_Y E$ to $J$ whose degree has only prime divisors in $S$. Thus,
  the level-$N$-structure $\alpha'$ will be mapped under $h'$ injectively into $J'\raro Y'$ so
  that $J'\raro Y'$ itself is equipped with a level-$N$-structure which we will call
  $\beta'$.\pagebreak[3]

  By Proposition (\ref{pro14}) and the choice of $Y'$, we see that $J'\raro Y'$ together with its
  principal polarization $\theta'$ (which is just a special kind of isogeny) are defined over
  $F'$. Also the level-$N$-structure $\beta'$ on $J'\raro Y'$ is defined over $F'$ since it is the
  image of the level-$N$-structure $\alpha'$ of $E'\times_{Y'}\cdots\times_{Y'} E'\raro Y'$ under
  $h'$. So, the triple $(J'\raro Y', \theta', \beta')$ is defined over $F'$.\pagebreak[3]

  The canonical morphism between fine moduli spaces
  $$ j^{(N)}: \cM^{(N)}_g \raro \cA^{(N)}_{g,1} $$
  which sends a curve with level-$N$-structure to its principally polarized Jacobian with the same
  level-$N$-structure is 2-to-1 over its image for $g\geq 3$ respectively 1-to-1 for $g=2$, see
  \cite{OortSteenbrink}. Hence, the curve $C'\raro Y'$ can be defined over $F'$ after applying a
  base change $Y''\raro Y'$ of degree at most 2.\pagebreak[3]

  As mentioned above, the choice of $Y'$ and $F'$ depend only on the family of elliptic curves
  $E\raro Y$.
\end{proof}\pagebreak[3]

So we see that a family of curves, whose Jacobian is isogenous to the $g$-fold product of a modular
family of curves, can be defined over a number field after a finite base change of degree at most
2. And the number field does not depend on the given family of curves but only on the modular family
of elliptic curves.\pagebreak[3]


\section{Reduction to finite fields}
 
Let $C\raro Y$ be a family of curves, i.\,e. a flat, projective morphism whose generic fiber is
smooth and geometrically connected. Assume that the Jacobian $J\raro Y$ is $Y$-isogenous to the
$g$-fold product $E\times_Y\cdots\times_Y E$ of a modular family of elliptic curves $E\raro
Y$.\pagebreak[3]

Because of Theorem (\ref{uniformDescend}) we may assume that our base field is a number field
$F$. We want to show that there is a finite prime of $F$ such that the family $\widetilde{C}\raro
\widetilde{Y}$  obtained by reduction modulo this prime has a also a smooth generic fiber. It is
clearly possible to find such a prime depending on the given family $C\raro Y$. But we want to show
that there is a choice of this prime which depends only on $E\raro Y$ and not on $C\raro Y$ or $g$,
so that for any family $C\raro Y$ of arbitrarily large genus $g$, whose Jacobian is $Y$-isogenous to
$E\times_Y\cdots\times_Y E$, its reduction $\widetilde{C}\raro\widetilde{Y}$ will still be
generically smooth.\pagebreak[3] 

Therefor we need a criterion for the generic smoothness of a family of curves $C\raro Y$ after
reduction modulo a prime. The characterization is given in terms of the existence of certain
endomorphisms on the \pp\ Jacobian $J\raro Y$ of $C\raro Y$.\pagebreak[3]

In this section we consider all schemes, morphisms and fiber products to live over an at most
$1$-dimensional base scheme $S$ which is suppressed from the notation.\pagebreak[3]

\begin{definition}[Split \ppav]
  A \ppav\ $(A,\lambda_A)$ {\it splits} if there are two positive-dimensional \ppavs\
  $(B,\lambda_B)$ and $(C,\lambda_C)$ such that $(A,\lambda_A)$ is isomorphic to $(B\times C,
  \lambda_B \times \lambda_C)$ as a \ppav.
\end{definition}\pagebreak[3]

For the next proposition, we assume that $S$ is the spectrum of an algebraically closed field of
arbitrary characteristic. 

Note that a smooth curve $C$ has a proper Jacobian $J$. But the converse is not true, since
e.\,g. two smooth curves intersecting transversally in one point also have a proper Jacobian. So we
need a criterion to distinguish between these cases.\pagebreak[3]

\begin{proposition}[Reducibility criterion for curves]\label{reducibilityCriterion}
  Let $C$ be a curve with proper Jacobian $(J,\lambda)$. Then $C$ is reducible if and only if
  $(J,\lambda)$ splits as a \ppav.
\end{proposition}\pagebreak[3]

\begin{proof}
  If $C$ has a proper Jacobian then it is either smooth or it consists of smooth irreducible
  components $C_i$ intersecting in a way such that they form a tree. Hence, if $C$ is reducible,
  then $(J,\lambda)$ is the product of the Jacobians $(J_i,\lambda_i)$ of the smooth components
  $C_i$.\pagebreak[3]

  It remains to show the converse that $C$ is reducible if $(J,\lambda)$ splits. Assume that
  $(J,\lambda) = (A_1\times A_2, \lambda_1\times\lambda_2)$ where $(A_i,\lambda_i)$ are
  positive-dimensional \ppavs\ and that $C$ is smooth. Choose an embedding
  $C\stackrel{f}{\hookrightarrow}J$ and let $p_i$ be the projection $A_1\times A_2 \rar{p_i}
  A_i$.\pagebreak[3]

  Then $C_i := p_{i*} f(C)\subset A_i$ is a $1$-cycle generating $A_i$ \big(i.\,e. $A_i$ is the
  smallest Abelian subvariety of $A_i$ containing $C_i$\big) because $C$ generates its Jacobian
  $J=A_1\times A_2$. Define $\tilde{C}:= C_1\times\{0\} + \{0\}\times C_2 \subset A_1\times
  A_2=J$. It follows that $\tilde{C}$ is a $1$-cycle generating $J$.\pagebreak[3]

  Let $\Theta_i\subset A_i$ be a divisor inducing the polarization
  $A_i\rar{\lambda_i}\wh{A_i}$. Then the divisor $\Theta := \Theta_1\times A_2 + A_1\times \Theta_2$
  on $A_1\times A_2=J$ induces the polarization $J\rar{\lambda}\wh{J}$. We compute the intersection
  number
  \begin{eqnarray*}
    \big(\tilde{C}.\Theta\big) &=& \phantom{+}\big(C_1\times\{0\}.\Theta_1\times A_2\big) +
                           \big(C_1\times\{0\}.A_1\times\Theta_2\big)\\
                       & & + \big(\{0\}\times C_2.\Theta_1\times A_2\big)
                           + \big(\{0\}\times C_2.A_1\times\Theta_2\big).
  \end{eqnarray*}
  The two middle terms are zero as an application of the projection formula shows. An other
  application of the projection formula on the remaining two terms gives us 
  $$ \big(\tilde{C}.\Theta\big) = \big(f(C).\Theta_1\times A_2\big) +
     \big(f(C).A_1\times\Theta_2\big) = \big(f(C).\Theta\big) = g $$ 
  where the last equality follows from the fact that $C$ is the Jacobian of $C$.\pagebreak[3]

  From the Matsusaka-Ran Theorem \cite{Ran-Matsusaka} follows that the $(A_i,\lambda_i)$
  are the Jacobians of the curves $C_i$ and that the $C_i$ are components of $C$, so that $C$ has to
  be reducible.
\end{proof}\pagebreak[3]

A characterization of split \ppav\ $(A,\lambda)$ is given in the next proposition. There, the
{\it Rosati-involution} on $\End_S(A)$ is denoted by
$$ f\mapsto f^\dag := \lambda^{-1}\circ\widehat{f}\circ\lambda $$
where $\widehat{A}\rar{\widehat{f}}\widehat{A}$ is the dual map of $f$.

\begin{proposition}[Splitting criterion]\label{splittingCriterion}
  For a \ppav\ $(A,\lambda_A)$, the following two statements are equivalent:
  \begin{ilist}
    \item $(A,\lambda_A)$ splits, i.\,e. $(A,\lambda_A)$ is isomorphic as a \ppav\ to a product
      $(B\times C,\lambda_B\times\lambda_C)$ of two positive-dimensional \ppavs\ $(B,\lambda_B)$ and
      $(C,\lambda_C)$.

    \item $(A,\lambda_A)$ possesses a non-trivial symmetric idempotent endomorphism, i.\,e. it
      exists a map $f\in\End_S(A)$ different from the identity and the zero map such that the two
      relations $f^\dag=f$ and $f^2=f$ hold.
  \end{ilist}
\end{proposition}\pagebreak[3]

\begin{proof}
  $(i)\Rightarrow(ii)$ Let $h:A\stackrel{\sim}{\longrightarrow}B\times C$ be an
  isomorphism of \ppavs\ and define $f$ to be the following composition of maps
  $$\xymatrix{
    {B\times C}\ar[r]^{1\times0}  &  {B\times C}\ar[d]^-{h^{-1}}\\
    {A}\ar[u]^-h\ar[r]^{=:f}  &  {A.}
  }$$
  Then $f$ is an idempotent and symmetric endomorphism of $A$. For the idempotence consider the
  commutative diagram
  $$\xymatrix{
    {B\times C}\ar[r]^{1\times0}\ar@/^2pc/[rr]^{1\times0} & {B\times
      C}\ar[r]^{1\times0}\ar@<-1ex>[d]_{h^{-1}} & {B\times C}\ar[d]^{h^{-1}} \\
    A\ar[u]^h\ar[r]^f\ar@/_1.5pc/[rr]_{f^2\stackrel{!}{=}f} & A\ar[r]^f\ar@<-1ex>[u]_{h} & A
  }$$
  where the lower row gives us $f^2$. If we follow the upper way around we get $f$. So $f$ and
  $f^2$ coincide and, therefore, $f$ is idempotent. For the symmetry look at the diagram
  $$\xymatrix{
    {\widehat{A}}\ar[r]^{\widehat{h}^{-1}}\ar@/^2.6pc/[rrr]^{\widehat{f}}  &
    {\widehat{B}\times\widehat{C}}\ar[r]^{\widehat{1}\times\widehat{0}}  &
    {\widehat{B}\times\widehat{C}}\ar[r]^{\widehat{h}}\ar[d]^{\lambda_B^{-1}\times\lambda_C^{-1}}  &  
    A\ar[d]^{\lambda_A^{-1}}  \\
    A\ar[u]^{\lambda_A}\ar[r]^-h\ar@/_1.5pc/[rrr]_{f\stackrel{!}{=}f^\dag}  &
    {B\times C}\ar[u]^{\lambda_B\times\lambda_C}\ar[r]^{1\times0}  &
    {B\times C}\ar[r]^-{h^{-1}}  &  A
  }$$
  which is commutative since dualizing endomorphisms commutes with inverting them. Again the
  lower row gives us $f$ while the upper way around we obtain $f^\dag$. So $f$ and $f^\dag$ are
  identical, telling us that $f$ is a symmetric and idempotent endomorphism of $A$.
  
  $(ii)\Rightarrow(i)$ Let $A\rar{f}A$ be a symmetric idempotent endomorphism of
  $A$. Define $B:=\Image(A\rar{f}A)$ and $C:=\Image(A\rar{1-f}A)$. Then we get a homomorphism
  $$ B\times C=fA\times(1-f)A \stackrel{h}{\longrightarrow} A,\quad (b,c)\mapsto b+c. $$
  Since $f$ is idempotent, the homomorphism
  $$ A\longrightarrow B\times C=fA\times(1-f)A,\quad a\mapsto \big(fa, (1-f)a\big) $$
  is an inverse map for $h$ and, therefore, $B\times C\stackrel{h}{\longrightarrow} A$ is an
  isomorphism of Abelian varieties.
  
  Let $\lambda_B$ and $\lambda_C$ be the restrictions of $\lambda_A$ on $B$ and $C$. This makes
  $B$ and $C$ into \ppavs. Since $f$ is symmetric the two diagrams
  $$\xymatrix{
    {\widehat{A}}\ar[r]^{\widehat{f}}  &  {\widehat{A}}  & {}\ar@{}[dr]|{\mathrm{ }} &  & 
    {\widehat{A}}\ar[r]^{\widehat{1}-\widehat{f}}  &  {\widehat{A}} \\ 
    A\ar[u]^{\lambda_A}\ar[r]^f  &  A\ar[u]_{\lambda_A}  &  &  &
    A\ar[u]^{\lambda_A}\ar[r]^{1-f}  &  A\ar[u]_{\lambda_A}
  }$$
  commute. In particular, $B$, which is the image of $f$, is mapped under $\lambda_A$ into the
  image of $\wh{A}$ under $\wh{f}$. The same holds for $C$ and $1-f$. But then also the diagram
  $$\xymatrix{
    {\widehat{B}\times\widehat{C}=\widehat{f}\widehat{A}\times(\widehat{1}-\widehat{f})\widehat{A}} & 
    {\widehat{A}}\ar[l]_-{\widehat{h}} \\
    {B\times C=fA\times(1-f)A}\ar@<-2ex>[u]^{\lambda_B\times\lambda_C}\ar[r]^-h  &  A\ar[u]_{\lambda_A}
  }$$
  commutes. Hence, as a \ppav\ $(A,\lambda_A)$ is isomorphic to $(B\times C,\lambda_B \times
  \lambda_C)$ via the map $h$.\
\end{proof}\pagebreak[3]

Assume that an Abelian scheme splits after reduction modulo a prime so that it owns a symmetric
idempotent endomorphism. To show that it is already split before reduction we want to lift the
endomorphism.\pagebreak[3]

Let $S = \Spec R$ be the spectrum of a henselian discrete valuation ring with quotient field $K$ and
residue field $k$. We use the following notational convention. A small subscript denotes the base
scheme. So the schemes $X_K$ and $X_k$ are schemes over $\Spec K$ and $\Spec k$, respectively. A
Scheme over $S$ is simply denoted by $X$ instead of $X_S$. Then $X_K$ is its general fiber and $X_k$
is its special fiber. Let $A\raro S$ be an Abelian scheme. The question we will study is when does
an endomorphism of $A_k$ lift to an endomorphism of $A_K$.\pagebreak[3]

\begin{definition}[The lifting property]
  We say that every endomorphism of $A_k$ {\it lifts} if the restriction map
  $$ 
  \begin{matrix}
    \End_S(A) & \longrightarrow & \End_k(A_k) \\
    f         & \longmapsto     & f_k := f_{|A_k}
  \end{matrix}
  $$
  is an isomorphism.
\end{definition}\pagebreak[3]

\begin{example}
  If $E\rar{}S$ is a relative elliptic curve such that $E_K$ and $E_k$ are both elliptic curves
  without complex multiplication so that the endomorphism rings $\End_S(E)$ and $\End_k(E_k)$ are
  isomorphic to $\bZ$, then the restriction map is an isomorphism because the endomorphisms
  are the multiplication-by-$m$ maps.\pagebreak[3] 

  The same is true for the $g$-fold product of the elliptic curve $E$ since in this case any
  endomorphism is build up from multiplication-by-$m$ maps which lift.
\end{example}\pagebreak[3]

The lifting property is invariant under \'etale isogenies.\pagebreak[3]

\begin{proposition}[The lifting property and \'etale isogenies]\label{liftingProperty}
  Let $A$ and $B$ be two \ppass\ over $S$ and $h:A\raro B$ an \'etale $S$-isogeny. If every
  endomorphism of $A_k$ lifts, then every endomorphism of $B_k$ lifts too.
\end{proposition}\pagebreak[3]

\begin{proof}
  Look at the commutative diagram
  $$
    \xymatrix{
      {\End_S(B)}\ar@{^{(}->}[r]^{h^*}\ar@{^{(}->}[d]  &  {\End_S(A)}\ar[d]^{\cong}    \\
      {\End_k(B_k)}\ar@{^{(}->}[r]^{h_k^*}               &  {\End_k(A_k)}
  }$$
  where $h^*$ is given by $f\mapsto h^*f := h^\dag\circ f\circ h$ and $h^*_k$ is the analogous map
  for endomorphisms of the special fiber. Let $f_k\in\End_k(B_k)$ be an endomorphism of $B_k$. We want
  to lift $f_k$ to an endomorphism $f\in\End_S(B)$ so that $f_{|B_k} = f_k$. Look at the map
  $h^*_kf_k\in\End_k(A_k)$. Since $A$ has the lifting property, the map $h^*_kf_k$ lifts to a map
  $A\rar{u}A$ so that $u_k = h^*_kf_k$. If we can show that $A\rar{u}A$ lies in the image of $h^*$,
  i.\,e. there is a map $f\in\End_S(B)$ with $h^*f = h^\dag\circ f\circ h = u$, then the map $f$
  is a lifting of $f_k$ because of the commutativity of the diagram above.\pagebreak[3]

  We know that $u_k = h^*_kf_k = h_k^\dag\circ f_k\circ h_k$ factorizes through $h_k$ so that
  $\Ker(h_k)$ is a subgroup scheme of $\Ker(u_k)$. Since our base $S$ is henselian and $\Ker(h_k)$
  \'etale, there is a subgroup scheme $G\subset\Ker(u)$ such that $G_k=\Ker(h_k)$. But then, being a
  subgroup scheme of $A$, the group scheme $G$ has to coincide with $\Ker(h)$ since finite \'etale
  schemes over $S$ are uniquely determined by their special fiber \cite[p.34]{MilneEC}. Hence, $u$
  factorizes through $h$, i.\,e. there is a map $B\rar{g}A$ such that $u = g\circ h$
  holds.\pagebreak[3] 

  Analogously one shows that the dual $\wh{g}$ of $g$ factorizes through the dual $\wh{h^\dag}$ of
  $h^\dag$. Hence, there exists an endomorphism $B\rar{f}B$ such that $g = h^\dag\circ f$ is
  valid. Therefore, we get the identity $u = h^*f$ and $f$ becomes a lifting of $f_k$. This implies
  that the Abelian scheme $B$ also has the lifting property.
\end{proof}\pagebreak[3]

In particular, Abelian schemes isogenous to a $g$-fold product of a relative elliptic curve without
complex multiplication have the lifting property.\pagebreak[3]

Let us return to our situation in the beginning of this section. Let $Y_F$ be a smooth, projective,
geometrically connected curve over some number field $F$. Let $\cO_F$ denote the ring of integers of
$F$ and let $\Spec \Fq\rar{}\Spec \cO_F$ be any finite point of $\Spec \cO_F$. We can extend
$Y_F\rar{}\Spec F$ to an arithmetic minimal model $Y\rar{}\Spec\cO_F$, i.\,e. $Y\rar{}\Spec\cO_F$ is
an integral, proper, regular, excellent and flat surface of finite type with general fiber $Y_F$
together with the usual minimality property similar to the geometric case, see
\cite{ChinburgMinimalModel}. Let $Y_{\Fq}\rar{}\Spec \Fq$ be the special fiber of $Y\rar{}\Spec
\cO_F$ over the point $\Spec \Fq\rar{}\Spec \cO_F$.\pagebreak[3]

Let $K$ be the function field of $Y_F$ so that $\Spec K\rar{}Y_F$ is the generic point. Let $k$ be
the function field of an irreducible component of $Y_{\Fq}$ so that $\Spec k\rar{}Y_\Fq$ is the
generic point of the corresponding irreducible component. Furthermore, let $R$ be the local ring of
$Y$ at this irreducible component. In particular, $R$ is a discrete valuation ring ($Y$ is regular)
with generic point $\Spec K\rar{}\Spec R$ and special point $\Spec k\rar{}\Spec R$.\pagebreak[3]
$$\xymatrix{
  {\{C_K,J_K,E_K\}\ }\ar@{^{(}->}[r]\ar[d]  &  {\{C,J,E\}}\ar[d]  &  
                                                     {\ \{C_k,J_k,E_k\}}\ar@{_{(}->}[l]\ar[d] \\ 
  {\Spec K\ }\ar@{^{(}->}[r]\ar[d] &  {\Spec R}\ar[d]    &  {\ \Spec k}\ar@{_{(}->}[l]\ar[d]   \\
  {Y_F\ }\ar@{^{(}->}[r]\ar[d]     &  {Y}\ar[d]          &  {\ Y_{\Fq}}\ar@{_{(}->}[l]\ar[d]    \\
  {\Spec F\ }\ar@{^{(}->}[r]       &  {\Spec \cO_F}      &  {\ \Spec \Fq}\ar@{_{(}->}[l]
}$$\label{diag}\pagebreak[3]

Finally, let $C_K$ be a smooth, projective, geometrically connected curve defined over the function
field $K$, $J_K$ its Jacobian and $E_K$ an elliptic curve. We may extend $C_K\rar{}\Spec K$ to a
minimal model $C\rar{}\Spec R$ with Jacobian $J\rar{}\Spec R$ and we denote the special fibers of
these models by $C_k\rar{}\Spec k$ and $J_k\rar{}\Spec k$. Of course, also $E_K\rar{}\Spec K$
extends to a (N\'eron) model $E\rar{}\Spec R$ with special fiber $E_k\rar{}\Spec k$.\pagebreak[3]

With a view towards Proposition (\ref{uniformDescend}), we want to allow our curve $C$ to be defined
over some finite field extension $K'$ of $K$ rather than over $K$ itself. Therefore, let $Y'\raro Y$
be the map of minimal arithmetic surfaces induced by $\Spec K'\raro\Spec K$, so that $K'$ is the
function field of $Y'$. Let $Y'_\Fq$ be the fiber of $Y'\raro \Spec\cO_F$ over the point $\Spec\Fq
\raro\Spec\cO_F$, $k'$ the function field of an irreducible component of $Y'_\Fq$ and $R'$ the local
ring of $Y'$ at this irreducible component. We can extend $C_{K'}$ to a model $C'\raro \Spec R'$ and
denote its special fiber by $C_{k'}$.\pagebreak[3]

\begin{theorem}[Reduction to finite fields]\label{reductionToFiniteFields}
  Let $E_K$ be a non-isotrivial elliptic curve. Then there is a finite point $\Spec \Fq \raro \Spec
  \cO_F$, depending only on $E_K$, such that the following property holds:

  Let $J_K$ be an Abelian
  variety which is $K$-isogenous to a $g$-fold product of $E_K$. Let $K'$ be a finite extension of
  $K$ such that $J_{K'}$ becomes the Jacobian of a projective and geometrically connected curve
  $C_{K'}$. Then $C_{K'}$ is smooth over $K'$ if and only if its reduction $C_{k'}$ with respect to
  $\Spec\Fq\raro\Spec\cO_F$ is smooth over $k'$.
\end{theorem}\pagebreak[3]

\begin{proof}
  First, we start with the case $K'=K$. We choose a finite point $\Spec \Fq\rar{}\Spec \cO_F$ of
  residue characteristic $p$ such that the following properties are fulfilled.
  \begin{1list}
    \item $Y_F$ has good reduction at $\Spec \Fq$, i.\,e. the fiber $Y_\Fq$ is a smooth curve. This
      depends only on $K$ -- the function field of $Y_F$ -- and is true for almost all points of
      $\Spec \cO_F$.

    \item $E_K\rar{}\Spec K$ extends to a smooth proper model $E\rar{}\Spec R$ such that
      $E_k\rar{}\Spec k$ is a non-isotrivial elliptic curve. This is true for almost all points of
      $\Spec \cO_F$ and depends only on $E_K$.

    \item There is an isogeny $E_K\times_K\cdots\times_K E_K \rar{}J_K$ such that its degree is
      prime to $p$. Using Proposition (\ref{pro15}) we see that this is true for almost all
      points of $\Spec \cO_F$ and depends only on $E_K$. Together with (2) this property will enable
      us to lift endomorphisms of $J_k$ to endomorphisms of $J_K$ with the help of Proposition
      (\ref{liftingProperty}).
  \end{1list}\vspace{1ex}\pagebreak[3]
  Since the three conditions above hold separately for all but finitely many points of $\Spec
  \cO_F$, we can find a point $\Spec \Fq\rar{}\Spec \cO_F$ fulfilling all conditions
  simultaneously. The choice of this point depends only on $E_K$.\pagebreak[3]

  As explained above, let $R$ be the local ring of $Y$ at $Y_\Fq$. Extend the curve
  $C_K\rar{}\Spec K$ to a minimal model $C\rar{}\Spec R$. Its Jacobian $J\rar{}\Spec R$ is equipped
  with a canonical principal polarization $\lambda: J\rar{}\wh{J}$ such that $(J_K,\lambda_K)$
  and $(J_k,\lambda_k)$ are the principally polarized Jacobian of $C_K$ and $C_k$,
  respectively. Since by assumption $J_K$ is $K$-isogenous to the $g$-fold product of $E_K$, the
  Jacobian $J_k$ is $k$-isogenous to the $g$-fold product of $E_k$ (actually $J\rar{}\Spec R$ is
  isogenous over $\Spec R$ to the $g$-fold product of $E\rar{}\Spec R$).\pagebreak[3]

  Let $\wh{R}$ be the completion of $R$ and $\wh{K}$ its quotient field. So after the base change
  $\Spec \wh{R} \rar{} \Spec R$, we get a model $C_{\wh{R}}\rar{}\Spec \wh{R}$ with generic fiber
  $C_{\wh{K}}\rar{}\Spec\wh{K}$ and special fiber $C_k\rar{}\Spec k$.\pagebreak[3]

  As discussed in the proof of Proposition (\ref{reducibilityCriterion}), a non-smooth curve with
  proper Jacobian becomes reducible after some finite base change. By Proposition
  (\ref{reducibilityCriterion}), the reducibility of the curve is equivalent to the splitting of its
  Jacobian. By Proposition (\ref{splittingCriterion}), the splitting is equivalent to the existence
  of a symmetric idempotent endomorphism of the Jacobian.\pagebreak[3]

  Since Proposition (\ref{liftingProperty}) tells us that $J_{\wh{K}}$ has the lifting property, a
  symmetric idempotent endomorphism exists on $J_{\wh{K}}$ if and only if it exists on $J_k$ (apply a
  finite base change $\Spec\wh{S}\rar{}\Spec\wh{R}$ if necessary). Hence $C_K$ is smooth over $K$
  if and only if $C_k$ is smooth over $k$.\pagebreak[3]

  Now, if $K'$ is a finite extension of $K$, with the same choice of prime
  $\Spec\Fq\raro\Spec\cO_F$, as above, the Jacobian $J_{\wh{K'}}$ still has the lifting property. So,
  replacing $K,R,k$ by $K',R',k'$, respectively, verbatim the same argument, as before, shows that
  $C_{K'}$ is smooth over $K'$ if and only if $C_{k'}$ is smooth over $k'$.
\end{proof}\pagebreak[3]

In particular, taking the minimal model $\widetilde{C}$ of $C_k$ over $\widetilde{Y} = Y_\Fq$, the
generically smooth family of curves $C\raro Y$ induces a corresponding generically smooth family of
curves $\widetilde{C}\raro \widetilde{Y}$ over a suitable finite field $\Fq$, which does not depend
on $C\raro Y$. In the next section, we see that the genus $g$ of the fibers of $\widetilde{C}\raro
\widetilde{Y}$, which is the same as the genus of the fibers of $C\raro Y$, is bounded.\pagebreak[3]


\section{Bounding the genus in characteristic p}

Let $\Fq$ be a finite field of characteristic $p$. Serre has shown in \cite{SerreRepartition} based
on the work of Tsfasman and Vl\u adu\c t \cite{TsVlAsympZeta} that for a curve $C/\Fq$, whose
Jacobian $J/\Fq$ is $\Fq$-isogenous to a product of elliptic curves, the genus $g$ is
bounded.\pagebreak[3]

We will compute an explicit bound for an easy special case and generalize the statement for families
of curves.\pagebreak[3]

\begin{proposition}[Explicit bound]\label{explicitBound}
  Let $C/\Fq$ be a curve of genus $g$ whose Jacobian $J/\Fq$ is $\Fq$-isogenous to the $g$-fold
  product of an elliptic curve $E/\Fq$. If $\TrF{E}{} > 0$, then $g\leq q+1$.
\end{proposition}\pagebreak[3]

\begin{proof}
  Isogenies become isomorphisms on $\ell$-adic cohomology, so that we have isomorphisms (where
  $\bar{C} = C\times_\Fq \bar{\bF}_q$, etc.)
  $$ \Hl{C}\cong\Hl{J}\cong H^1(\bar{E}\times\cdots\times\bar{E}, \Q_\ell) \cong 
     \bigoplus_{i=1}^g\Hl{E} $$
  compatible with the action of Galois. Therefore, the Weil conjectures tell us that the number of
  $\Fq$-rational points of $C$ is given by
  \begin{eqnarray*}
    0\leq \#C(\Fq) &=& q+1 - \TrF{C}{} \\
                   &=& q+1 - g\cdot\TrF{E}{}
  \end{eqnarray*}
  where $F_C$ and $F_E$ denote the $q$-th power Frobenius on $C$ and $E$, respectively. Hence, it
  follows that
  $$ g \leq \frac{q+1}{\TrF{E}{}} \leq q+1 $$
  giving the desired bound.
\end{proof}\pagebreak[3]

With the same method, one can compute explicit bounds for other Frobenius traces, too, but we don't
need to.\pagebreak[3]

Now, let $Y$ be a smooth, projective, geometrically connected curve over some  finite field $\Fq$
and $C\rar{}Y$ a semistable family of curves of genus $g$. To bound the genus $g$, we will count the
minimal number $\delta$ of singularities in the geometric fibers of $C\raro Y$ and compare it with
the following natural upper bound.\pagebreak[3]

\begin{proposition}[Upper bound for $\delta$]\label{upperBound}
  Let $f:C\raro Y$ be a semistable family of curves of genus $g\geq 2$. Then
  $$ \delta\leq 12\cdot\deg f_*\omegacy $$
  where $\delta$ is the number of singularities in the geometric fibers of $f$ and $\omegacy$ is the
  relative dualizing sheaf.
\end{proposition}\pagebreak[3]

\begin{proof}
  The characteristic $p$ case was proven by Szpiro in \cite[Prop.1]{Szpiro}.
\end{proof}\pagebreak[3]

If the Jacobian $J\raro Y$ of $C\raro Y$ is $Y$-isogenous to the $g$-fold product of a
non-isotrivial family of elliptic curves $E\rar{}Y$, then we can express this upper bound in terms
of the genus $g$ and the height of $E\raro Y$.\pagebreak[3]

Recall that the {\it height} of a group scheme $G\raro Y$ is given by
$$ h(G) := \deg(s^*\Omega^1_{G/Y}) $$
where $s:Y\raro G$ is the zero section.\pagebreak[3]

\begin{corollary}[Bounding $\delta$ for decomposable Jacobians]\label{upperBound2}
  Let $C\raro Y$ be a semistable family of curves whose Jacobian $J\rar{}Y$ is $Y$-isogenous to the
  $g$-fold product of a non-isotrivial family of elliptic curves $E\rar{}Y$. Then
  $$ \delta \leq 12\,h(E)\cdot g $$
  where $h(E)$ is the height of $E\rar{}Y$. In particular, the constant $h(E)$
  depends only on $E\rar{}Y$, so that $\delta$ is linearly bounded by $g$.
\end{corollary}\pagebreak[3]

\begin{proof}
  Let $E\times_Y\cdots\times_Y E\raro J$ be an isogeny with kernel $N$. Then we have $h(J) =
  h(E\times_Y\cdots\times_Y E) - h(N)$. Since $E$ is non-isotrivial, $N$ is an extension of an
  \'etale group scheme by some factors of the form $\mu_{p^n}$. Both group schemes have height zero
  (the Cartier dual of $\mu_{p^n}$ is \'etale), so that $h(N)=0$.\pagebreak[3]

  Furthermore, we have $h(E\times_Y\cdots\times_Y E) = g\cdot h(E)$ since the sheaf
  $\Omega^1_{E\times\cdots\times E/Y}$ is isomorphic to $\bigoplus_{i=1}^g p_i^*\Omega^1_{E/Y}$
  where $p_i$ is the projection on the $i$-th factor.\pagebreak[3]

  The height of the Jacobian is related to $C\raro Y$ by $h(J)=\deg f_*\omegacy$ because
  $\det s^*\Omega^1_{J/Y}\cong\det f_*\omegacy$ \cite[p.351]{Faltings}. So we have
  $$ \delta\leq 12\cdot\deg f_*\omegacy=12\cdot h(J) = 12h(E)\cdot g $$
  where the inequality is given by Proposition (\ref{upperBound}).
\end{proof}\pagebreak[3]

This gives us an upper bound for $\delta$ in terms of $g$. For the lower bound we will use
Proposition (\ref{explicitBound}). To determine how many fibers of $E\raro Y$ have a positive
Frobenius trace, we will use the Sato-Tate conjecture.\pagebreak[3]

\begin{theorem}[Sato-Tate Conjecture]\label{SatoTate}
  Let $E\rar{}Y$ be a non-isotrivial family of elliptic curves and $a$ and $b$ two real numbers
  between $0$ and $\pi$. Then
  $$ \lim_{n\rightarrow\infty} \frac{\#\big\{y\in Y(\Fqn)\ \big|\ a\leq\Theta(y)\leq b\big\}}{q^n}
     = \frac{2}{\pi} \int_a^b \sin^2\varphi\, d\varphi $$
  where $\Theta(y)$ is the angle of a Frobenius eigenvalue of the fiber $E_y$, i.\,e. the
  eigenvalues of the Frobenius acting on $\Hl{E_y}$ are given by $q^{n/2}\cdot e^{\pm\Theta(y)i}$.
\end{theorem}\pagebreak[3]

\begin{proof}
  This was proven by Deligne in \cite[p.212, (3.5.7)]{Weil2}.
\end{proof}\pagebreak[3]

We derive the following lower bound for the number of singularities $\delta$. Since later we want to
apply this Proposition together with Proposition (\ref{uniformDescend}) we allow the families of
curves to be defined over some finite covering $Y'\raro Y$.\pagebreak[3]

\begin{proposition}[Lower bound for $\delta$]\label{lowerBound}
  Let $E\raro Y$ be a non-isotrivial family of elliptic curves and $Y'\raro Y$ some finite
  covering of degree $d$. Let $C'\raro Y'$ be a semistable family of curves of genus $g$ whose
  Jacobian $J'\raro Y'$ is $Y'$-isogenous to the $g$-fold product of $E'\raro Y'$. Then there is a
  constant $c = c(E/Y, d) > 0$, depending only on $E\raro Y$ and $d$ such that
  $$ c\cdot \frac{\log g}{\log\log g}\cdot g \leq \delta $$
  where $\delta$ is the number of singularities in the geometric fibers of $C'\raro Y'$. In
  particular, $\delta$ is not linearly bounded above by $g$ and the lower bound does not depend on
  the particular choice of the covering $Y'\raro Y$.
\end{proposition}\pagebreak[3]

\begin{proof}
  We consider the case $d=1$ first. After enlarging $q$ if necessary, using the
  Sato-Tate-conjecture, we may assume that 
  $$\#\Big\{y\in Y(\Fqn)\, \Big|\, 0\leq \Theta(y)  <\frac{\pi}{2} \Big\} > \frac{1}{4}q^n. $$ 
  If $g>q^n+1$, then a fiber over an $\Fqn$-rational point
  $y$ of $Y$ with $\Theta(y)<\frac{\pi}{2}$ has to be singular by (\ref{explicitBound}). Its
  Jacobian is either isogenous to the $g$-fold product of a single elliptic curve or a torus. In the
  toric case, the curve has at least $g$ singularities. In the compact case, the curve is a chain of
  smooth curves each of genus less or equal to $q^n+1$. Such a curve will have at least
  $\big\lfloor\frac{g}{q^n+2}\big\rfloor$ singularities. Underestimating the number of
  singularities, we can say that in any case we have at least $\frac{g}{2q^n}$ singularities. So the
  total number of singularities we get from these fibers is at least 
  $$ \frac{1}{4}q^n\cdot\frac{g}{2q^n} = \frac{1}{8}g $$
  singularities.\pagebreak[3]

  There is one point we have to take care of. If $m$ is a natural number dividing $n$, then
  $Y(\F_{q^m}) \subset Y(\Fqn)$. So saying that we get $\frac{1}{8}g$ singularities from the
  $\F_{q^m}$-rational fibers and additional $\frac{1}{8}g$ singularities from the $\Fqn$-rational
  fibers is not fully correct because we possibly count some singularities more than once. To deal
  with this problem we will only consider extensions $\F_{q^e}$ of prime degree $e$.\pagebreak[3]

  Hence assume that $g-1 > q^2, q^3, q^5, q^7, q^{11}, \ldots, q^e,\ldots$ where the exponents $e$ are
  prime numbers. How many $q^e<g-1$ with $e$ prime are there? It is the number of primes $e$ with $e
  < \log_q(g-1)$. So by the prime number theorem, there is a constant $c_1>0$ such that there are
  at least $c_1\frac{\log_q(g-1)}{\log\log_q(g-1)}$ such primes $e$ ($c_1$ is a little bit less
  than $1$ if $g$ is large). So we get not less than
  $$ c_1\frac{\log_q(g-1)}{\log\log_q(g-1)} \cdot \frac{1}{8}g $$
  singularities up to multiply counted ones.\pagebreak[3]

  Thus, we have to deal with the singularities we counted more than once, namely the ones coming
  from fibers defined over $\Fq$-rational points because $Y(\Fq)\subset Y(\Fqn)$ for all $n$. Using
  a bad estimate for $\# Y(\Fq)$, we assume that there are at most $2q$ $\Fq$-rational points ($q$
  not too small, enlarge if necessary). Then we counted at most $2q\cdot\frac{g}{2q^e} =
  \frac{g}{q^{e-1}}$ points too often for each prime $e$. So an upper bound for the total error is
  $$ \sum_{e\ \mathrm{prime}} \frac{g}{q^{e-1}} \leq 2g. $$
  Therefore, the corrected total number of singularities we counted is
  $$ \left( \frac{c_1}{8} \frac{\log_q(g-1)}{\log\log_q(g-1)}-2 \right)g. $$
  Remember that we enlarged $q$ to apply the Sato-Tate-conjecture for $E\rar{}Y$. So for our
  original $q$, we can say that there is a constant $c>0$ depending only on $E\rar{}Y$ such that
  there are at least $c \frac{\log g}{\log\log g} g$ singularities.

  Now assume that we have a covering $Y'\raro Y$ of degree $d\geq 1$. Then any $\Fqn$-rational point
  of $Y$ has at least one $\F_{q^{rn}}$-rational preimage with $r\leq d$. So applying the 
  Sato-Tate-conjecture on $E'\raro Y'$, which is the extension of $E\raro Y$ with respect to the
  base change $Y'\raro Y$, we see that
  $$ \#\left\{y\in Y'(\F_{q^{dn}})\, \Big|\, 0\leq \Theta(y) <\frac{\pi}{2} \right\} \geq
     \#\left\{y\in Y(\Fqn)\, \Big|\, 0\leq \Theta(y) <\frac{\pi}{2d} \right\} > \varepsilon q^n $$
  where $\varepsilon >0$ is some constant depending only on $d$. Now verbatim the same counting as
  above gives us a constant $c=c(E/Y,d)$, depending only on $E\raro Y$ and $d$, such that
  $C'\raro Y'$ has at least $c\cdot\frac{\log g}{\log\log g}$ singularities.
\end{proof}\pagebreak[3]

It follows that the genus $g$ of such families of curves has to be bounded.\pagebreak[3]

\begin{theorem}[The genus is bounded]\label{genusBound}
  Let $E\raro Y$ be a non-isotrivial family of elliptic curves and $C'\raro Y'$ a family of curves
  of genus $g$ defined over a covering $Y'\raro Y$ of degree at most $d$. Assume that the Jacobian
  $J'\raro Y'$ of $C'\raro Y'$ is $Y'$-isogenous to the $g$-fold product of $E'\raro Y'$. Then the
  genus of $C'\raro Y'$ is bounded, i.\,e. there is a constant $B=B(E/Y, d)$, depending only on
  $E\raro Y$ and $d$, such that $g$ is smaller than $B$.
\end{theorem}\pagebreak[3]

\begin{proof}
  \Wlog, we may assume that $E\raro Y$ has semistable reduction everywhere. If not, we can achieve
  this after a finite base change using the semistable reduction theorem.

  Let $J'\raro Y'$ be the family of Jacobians of $C'\raro Y'$. By assumption $J'\raro Y'$ is
  $Y'$-isogenous to the $g$-fold product of $E'\raro Y'$. Hence, $J'\raro Y'$ has semistable
  reduction and, therefore, the family of curves $C'\raro Y'$ is semistable. So by
  (\ref{lowerBound}) and (\ref{upperBound}) the number $\delta$ of singularities in the geometric
  fibers of $C\raro Y$ satisfies (notice that $h(E')\leq d\cdot h(E)$)
  $$ c_0 \frac{\log g}{\log\log g}\cdot g \leq \delta \leq 12 d\cdot h(E)\cdot g $$ 
  where the constant $c_0>0$ depends only on $E\raro Y$ and $d$. But then $g$ cannot be arbitrarily
  large, since the left hand side is not linearly bounded by $g$. So there is a constant $B$
  depending only on $E\raro Y$ and $d$ such that $g$ is smaller than $B$.
\end{proof}\pagebreak[3]


\section{Conclusion}

We come to the results announced in the introduction. Our base field is now $\bC$
again.\pagebreak[3]

\begin{theorem}[Bound for the genus]\label{theorem2neu}
  Let $C\raro Y$ be a family of curves of genus $g$ whose Jacobian $J\raro Y$ is $Y$-isogenous to
  the $g$-fold product of a non-isotrivial family of elliptic curves $E\raro Y$ which can be defined
  over a number field. Then the genus $g$ is bounded, i.\,e. there is a number $B=B(E/Y)$
  depending only on $E\raro Y$ such that $g$ is smaller than $B$.
\end{theorem}\pagebreak[3]

\begin{proof}
  Using Proposition (\ref{uniformDescend}) we see that after replacing $Y$ by a finite covering,
  there is a covering $Y'\raro Y$ of degree at most $2$ such that $C'\raro Y'$ can be defined over a
  number field $F$ which depend only on $E\raro Y$.\pagebreak[3]

  Theorem (\ref{reductionToFiniteFields}) tells us that there is a finite prime $\Spec\Fq\raro
  \Spec\cO_F$, whose choice depends only on $E\raro Y$, such that the reduction of $C'\raro Y'$
  modulo this prime yields a family $\widetilde{C'}\raro\widetilde{Y'}$ of curves of genus $g$ over
  $\Fq$ whose Jacobian is $\widetilde{Y'}$-isogenous to the $g$-fold product of the reduction
  $\widetilde{E'}\raro\widetilde{Y'}$ of $E\raro Y$.\pagebreak[3]

  By Theorem (\ref{genusBound}) applied to $\widetilde{C'}\raro\widetilde{Y'}$ there is a number
  $B$, depending only on $E\raro Y$, such that $g$ is smaller than $B$.
\end{proof}\pagebreak[3]

Coming to Shimura curves, we first give a finiteness statement about modular families of elliptic
curves.\pagebreak[3]

\begin{proposition}[Finiteness of modular families of elliptic curves]\label{modularFiniteness}
  Fix two integers $q$ and $s$. Then there are only finitely many semistable modular families of
  elliptic curves $E\raro Y$ defined over a base curve $Y$ of genus at most $q$ and smooth outside a
  set $S\subset Y$ of cardinality at most $s$.
\end{proposition}\pagebreak[3]

\begin{proof}
  Let $E\raro Y$ be a modular family as in the proposition and let $Y\grar{j_E}\Pone$ be the $j$-map
  corresponding to the family $E\raro Y$. Because of the semistability of $E\raro Y$, an application
  of the ABC-conjecture for function fields yields 
  $$ \deg(j_E) \leq 6\cdot (2q-2+s) =:d. $$
  So, in particular, the degree of the $j$-map is absolutely bounded by $d$. Therefore, the modular
  family of elliptic curves $E\raro Y$ is given by a subgroup $\Gamma\subset\SL_2(\Z)$ of index at
  most $d$.\pagebreak[3]

  Since $\SL_2(\Z)$ is finitely generated, there are only finitely many subgroups $\Gamma$ of
  $\SL_2(\Z)$ of index at most $d$. Thus, we have only finitely many semistable modular families of
  elliptic curves $E\raro Y$ over a curve of genus at most $q$ and smooth outside a set of
  cardinality at most $s$. 
\end{proof}\pagebreak[3]

\begin{example}
  Let $q=0$ so that $Y=\Pone$. Beauville showed that for $s\leq3$ there are no
  non-isotrivial semistable families of elliptic curves at all. And for $s=4$, Beauville showed
  that there are exactly six non-isotrivial semistable families of elliptic curves, all
  modular, corresponding to the congruence subgroups $\Gamma(3)$, $\Gamma_1(4)\cap\Gamma(2)$,
  $\Gamma_1(5)$, $\Gamma_1(6)$, $\Gamma_0(8)\cap\Gamma_1(4)$ and $\Gamma_0(9)\cap\Gamma_1(3)$, see
  \cite{BeauvilleP1}.
\end{example}\pagebreak[3]

So, if the family of Jacobians is $Y$-isogenous to the $g$-fold product of a modular family of
elliptic curves, then we get the following result from Theorem (\ref{theorem2neu}).\pagebreak[3]

\begin{corollary}[Uniform bound for modular families]\label{cor1}
  Fix two integers $q$ and $s$. Then there is a constant $B=B(q,s)$ such that for any
  semistable family of curves $C\raro Y$, which is defined over a base curve $Y$ of genus at most
  $q$ and whose family of Jacobians $J\raro Y$ is smooth outside a set $S\subset Y$ of cardinality
  at most $s$ and $Y$-isogenous to the $g$-fold product of a modular family of elliptic curves
  $E\raro Y$, the genus $g$ of the fibers of $C\raro Y$ is bounded above by $B$. In particular,
  $B$ depends only on $q$ and $s$.
\end{corollary}\pagebreak[3]

\begin{proof}
  By Proposition (\ref{modularFiniteness}), there are only finitely many semistable modular families
  of elliptic curves $E\raro Y$ over a curve of genus at most $q$ and smooth outside a set of
  cardinality at most $s$. Because of the modularity, each one can be defined over some number field
  \cite{Corvallis}.\pagebreak[3]

  So Theorem (\ref{theorem2neu}) gives for each $E\raro Y$, a bound $B(E/Y)$ such that for $g$
  larger than $B(E/Y)$, the $g$-fold product of $E\raro Y$ is not $Y$-isogenous to a Jacobian. Thus,
  taking $B=B(q,s)$ to be the maximum of these finitely many numbers $B(E/Y)$ proves the corollary.
\end{proof}\pagebreak[3]

We derive the following corollary.\pagebreak[3]

\begin{corollary}[Curves over $\bP^1_\C$ with strict maximal Higgs field]
  Given an integer $s\geq 0$, there is a natural number $B=B(s)$, depending only on $s$, such
  that for any semistable family of curves $C\raro\bP^1_\C$, whose Jacobian $J\raro\bP^1_\C$ is
  smooth outside a set $S\subset\Pone$ of cardinality at most $s$ and has a strictly maximal Higgs
  field, the genus $g$ of the fibers of $C\raro\Pone$ is bounded above by the number $B$.
\end{corollary}\pagebreak[3]

\begin{proof}
  Because of Corollary (\ref{corollarySplit}), we may choose $B=B(s)$ to be the constant $B(0,s)$
  from Corollary (\ref{cor1}).
\end{proof}\pagebreak[3]

In particular, a given rational Shimura curve parameterizing a family of high-dimen\-sional Abelian
varieties with strictly maximal Higgs field does not lie in the closure of the Schottky
locus. This proves Theorem (\ref{theorem1}).\pagebreak[3]


\end{document}